\documentclass[11pt]{amsart}
\usepackage{amsmath, amssymb, amscd, amsthm, euscript}

\newtheorem{thm}{Theorem}[section]
\newtheorem{cor}[thm]{Corollary}

\newtheorem{lem}[thm]{Lemma}
\theoremstyle{definition}
\newtheorem{defn}[thm]{Definition}

\newtheorem{rem}[thm]{Remark}

\title[Maximal Conv. Groups and Rank One Symm. Spaces]{Maximal Convergence Groups and Rank One Symmetric Spaces}
\author{Ara Basmajian}
\address{Department of Mathematics, University of Oklahoma, Norman,
OK 73019} \email{abasmajian@ou.edu}

\author{Mahmoud Zeinalian}
\address{Department of Mathematics, C.W. Post Campus, Long Island
University, 720 Northern Boulevard, Brookville, NY 11548}
\email{mzeinalian@liu.edu}

\begin{document}
\bibliographystyle{h-elsevier2}
\maketitle

\begin{abstract}
We show  that the group of conformal homeomorphisms of the
boundary of a rank one symmetric space (except the hyperbolic
plane) of noncompact type acts as a maximal convergence group.
Moreover, we show that any family of uniformly quasiconformal
homeomorphisms has the convergence property. Our  theorems
generalize results of Gehring and Martin in the real hyperbolic
case for M\"obius groups. As a consequence, this shows that the
maximal convergence subgroups of the group of self homeomorphisms
of the $d$-sphere are not unique up to conjugacy. Finally, we
discuss some implications of maximality.

\end{abstract}

\section{Introduction}

The convergence property (Definition
\ref{defn:convergenceproperty}) is an essential property that all
families of M\"obius transformations possess \cite{GM}. Many of
the basic theorems in the theory of Kleinian groups can be proven
within this topological context. Quasiconformal and convergence
families have been studied in various contexts, including the
articles, \cite{GM} and \cite{T}. In \cite{GM}, Gehring and Martin
showed that, for $d\geq 2$, the M\"obius group acting on the
$d$-sphere is maximal in the sense that it is a convergence group
and that no group of homeomorphisms of the boundary sphere having
the convergence property can properly contain it. Since the
isometry groups of rank one symmetric spaces of noncompact type
are natural generalizations of the M\"obius groups, it is of
interest to see if these groups possess the same convergence and
maximality property.
We show the following,\\

\setcounter{section}{3} \setcounter{thm}{3}

\begin{thm} Let $\mathbb{H}$ be a rank one symmetric space of noncompact
type which is not the real hyperbolic plane. Then the  group of
conformal homeomorphisms of $\partial\mathbb{H}$ is a convergence
group. In fact, any family $\mathcal F$ of $K$-quasiconformal
homeomorphisms of $\partial\mathbb{H}$ has the convergence
property.
\end{thm}

For real hyperbolic space of dimension greater than two, the above
theorem was proven by Gehring and Martin in  \cite{GM}. Our proof
of the general case (for all rank one symmetric spaces of
noncompact type) applied to real hyperbolic space is significantly
different than theirs. See Tukia \cite{T} for related results.\\

\begin{thm} Let $\mathbb{H}_{\mathbb C}$ denote complex hyperbolic space
of complex dimension greater than one. Then the group of conformal
homeomorphism of its boundary, $\text{Conf}(\partial
\mathbb{H}_{\mathbb C})$, is a maximal convergence group. That is,

\begin{enumerate}
\item  $\text{Conf}(\partial \mathbb{H}_{\mathbb C})$ is a
convergence group, and

\item If $G \supseteq \text{Conf}(\partial \mathbb{H}_{\mathbb
C})$ is a convergence group acting on $\partial
\mathbb{H}_{\mathbb C}$, then $G=\text{Conf}(\partial
\mathbb{H}_{\mathbb C})$.
\end{enumerate}
\end{thm}

Putting theorem (\ref{thm:maximal}) together with the results of
Gehring, Martin \cite{GM} and Pansu \cite{P3}, as a corollary we
have,\\

\begin{cor} Let $ \mathbb H$  be  a rank one symmetric space of
noncompact type  different from the real hyperbolic plane. Then
the group of conformal homeomorphisms, $\text{Conf}(\partial
\mathbb{H})$, is a maximal convergence group.
\end{cor}

Note that this shows that the maximal convergence subgroups of the
group of self homeomorphisms of the $d$-sphere, $Homeo(S^d)$, are
not unique up to conjugacy. In fact, for infinitely many values of
$d$, the topological group $Homeo(S^d)$ contains maximal
convergence subgroups of different dimensions. For example, the
ten dimensional $\text{Isom}(\mathbb{H}_\mathbb{R}^4)$ and the
eight dimensional $\text{Isom}(\mathbb{H}_\mathbb{C}^4)$ act on
the $3$-sphere  as maximal convergence groups.

In  section four, we mention some elementary consequences of
corollary (\ref{cor:maximal}) to convergence groups and
quasi-isometry groups  acting on $\mathbb{H}\cup
\partial\mathbb{H}$.

The reason for  excluding the hyperbolic plane from the above
discussion is that unlike the higher dimensional situation, the
hyperbolic plane has trivial conformal structure on the boundary.
For instance, all smooth diffeomorphisms of the circle are
conformal. In this case, the affine structure of the boundary and
quasisymmetric mappings may be utilized, see \cite{BZ}.

Finally, in the literature some authors use the term convergence
group to mean a discrete convergence group. It is important to
note that such a restriction is not imposed in this paper.

\setcounter{section}{1}

\section{Basics}

Cartan's  classification of semisimple Lie groups implies that
there are three families of rank one symmetric spaces of
noncompact type and an exceptional one: real hyperbolic space
$\mathbb{H}_\mathbb{R}^n$, complex hyperbolic spaces
$\mathbb{H}_\mathbb{C}^{2n}$, quaternionic hyperbolic space
$\mathbb{H}_\mathbb{H}^{4n}$, and the Cayley plane
$\mathbb{H}_{\mathbb{C}a}^{16}$ (see \cite{Mos}). Throughout this
paper $\mathbb{H}$ stands for a rank one symmetric space of
noncompact type except the hyperbolic plane (or equivalently
complex hyperbolic space of complex dimension one).
 The boundary of $\mathbb{H}$, denoted by $\partial
\mathbb{H}$, is a smooth sphere which is naturally endowed with a
conformal structure. $\text{Isom}(\mathbb{H})$ denotes the
isometry group of $\mathbb{H}$ and
$\text{Conf}(\partial\mathbb{H})$ stands for  the group of
conformal homeomorphisms of the ideal boundary. Isometries of
$\mathbb{H}$ induce conformal homeomorphisms on $\partial\mathbb
H$. That is, the image of the natural map $\phi:
\text{Isom}(\mathbb{H})\rightarrow\text{Homeo}(\partial
\mathbb{H})$ defined by $f \mapsto f|_{\partial\mathbb{H}}$ is
precisely $\text{Conf}(\partial \mathbb{H})$. In fact, $\phi:
\text{Isom}(\mathbb{H})\rightarrow \text{Conf}(\partial
\mathbb{H})$ is an isomorphism. The stabilizer of a point $x \in
\mathbb H$ is a maximal compact subgroup of
$\text{Isom}(\mathbb{H})$ which we denote by $\mathcal{K}$. The
Iwasawa decomposition says $\text{Isom}(\mathbb{H})
=\mathcal{K}\mathcal{A}\mathcal{N}$; where $\mathcal{A}$ is a
one-parameter group of translations along the geodesic connecting
$x$ to $\infty$, and $\mathcal{N}$ is a nilpotent group (see
Section (29) of \cite{P2}). Each such $\mathcal{N}$ has a fixed
point $\infty \in
\partial \mathbb{H}$, and acts simply transitively on $\partial
\mathbb{H}-\{\infty \}$. Thus, the ideal boundary,
$\partial\mathbb H$, can be identified with the one point
compactification of $\mathcal{N}$, where $\mathcal{N}$ has a
naturally defined (Euclidean metric in the real case)
Carnot-Caratheodory metric, denoted $d_c$,  on it. With respect to
this metric, the elements of $\mathcal{N}$ act
  as isometries, and the elements of $\mathcal{A}$ as
dilations. Thus  associated to each $g \in \mathcal{A}$  is a
positive real number (its dilation factor) $\lambda (g)$
satisfying, $d_c(g(x),g(y))=\lambda(g) d_c(x,y)$, for all $x,y \in
\mathcal{N}$. Finally we note that the group of conformal
homeomorphisms of $\partial\mathbb H$ acts {\it almost triply
transitively} on points. That is, given three distinct points
$x,y, z \in
\partial\mathbb H$, and a positive real number $r$, there exists
an element of $\text{Conf}(\partial \mathbb{H})$ which takes $x$
to the identity element $e \in \mathcal{N}$, $y$ to $\infty$, and
$z$ to the Carnot-Caratheodory sphere of radius $r$ centered at $e
\in \mathcal{N}$ or any other fixed sphere. This fact, follows
from the above description of the actions of $\mathcal{N}$ and
$\mathcal{A}$ on $\mathbb{H}$. For a discussion on some of the
basics of the complex case, see \cite{BM} and \cite{KR}.

Throughout this paper all sequences are assumed to be infinite. As
a matter of convention, a subsequence is often denoted with the
same notation as the original sequence. We use $\rho(\cdot,\cdot)$
to denote the spherical distance on $S^d$.

\begin{defn}\label{defn:convergenceproperty}  Let $Y$ be a compact topological space.
A family $\mathcal F \subset \text{Homeo}(Y)$  is said to have the
{\it convergence property} if each infinite sequence $\{f_n\}$ of
$\mathcal F$ contains a subsequence which

\begin{description}
\item[(C1)] converges uniformly to an element of
$\text{Homeo}(Y)$, or

\item[(C2)] has the attractor-repeller property, that is, there
exists a point $a \in Y$, the {\it attractor}, and a point $r \in
Y$, the {\it repeller}, so that the $\{f_n\}$ converge to the
constant function $a$, uniformly outside of any open neighborhood
of $r$. Note that $a$ may equal $r$.
\end{description}

\end{defn}

\begin{rem}\label{rem:inverses for convergence property}
  If the sequence $\{f_n \}$ uniformly converges
to $f$, then $\{f_n^{-1}\}$ uniformly converges to $f^{-1}$.
Similarly, if $\{f_n\}$ has attractor $a$ and repeller $r$, then
$\{f_n^{-1}\}$ has attractor $r$ and repeller $a$.
\end{rem}

Discreteness for a convergence group is equivalent to    any
infinite sequence having  a subsequence for which axiom (C2)
holds.

\begin{defn} Let $(X, d)$ be a metric measure space.  For a map $f: X
\rightarrow X$, $x\in X$, and $r>0$ define

$$L_f(x,r)=\text{Max}\{d(f(x),f(y))|d(x,y)=r\}$$
$$l_f(x,r)=\text{Min}\{d(f(x),f(y))|d(x,y)=r\}$$
$$H_f(x, r)= \frac{L_f(x,r)}{l_f(x,r)}$$
$$H_f(x)= \limsup_{r\rightarrow 0} H_f(x, r)$$

A homeomorphism $f$ is called {\it quasiconformal}, if $H_f(x)$ is
uniformly bounded. In the presence of a measure, we say $f$ is
$K$-{\it quasiconformal} , if it is quasiconformal and for almost
every $x \in X$, $H_f(x) \leq K$. The map $f$ is said to be
$K'$-{\it quasisymmetric} if $H_f(x, r) \leq K'$ for all $x \in X$
and all $r >0$. A family is called {\it uniformly} $K$-{\it
quasiconformal}
 or {\it uniformly} $K^{\prime}$-{\it quasisymmetric} if there are constants $K$ or $K'$ for which all maps in the family are $K$-
 quasiconformal or $K'$-quasisymmetric, respectively.
\end{defn}

The above definition of a quasiconformal map is usually referred
to as the metric definition.  It follows from Theorem 10.19 of
\cite{H}, that a homeomorphism $f: \mathcal{N} \rightarrow
\mathcal{N}$ is $K'$-quasisymmetric if and only if there exists a
homeomorphism $\eta: [0, \infty) \rightarrow [0, \infty)$ such
that $$d(x, a) \leq td(x, b) \phantom{p} \text{implies}
\phantom{h} d(f(x), f(a)) \leq \eta(t) . d(f(x), f(b))$$ for every
$a, b, x \in X$ and for every $t\in [0, \infty)$. The
homeomorphism $\eta$ and $K'$ depend only on each other and data
associated to the space.

Note that a $K$-quasisymmetric map is always $K'$-quasiconformal,
for $K'=K$. Conversely,  for a large class of metric spaces
including the boundaries  of rank one symmetric spaces of
noncompact type, a $K$-quasiconformal map is $K'$-quasisymmetric,
where $K'$ depends only on $K$ and the data associated to the
ideal boundary. See Lemma 4.6 and Lemma 4.8 of \cite{HK} for
details. This implies that on these metric spaces a uniformly
$K$-quasiconformal family is also uniformly
$K^{\prime}$-quasisymmetric.

\section{Maximality of Conformal Groups}

The following lemma will allow us to make several simplifications,
through normalization, in our arguments involving the convergence
property.

\begin{lem}\label{lem:normalization1}
Suppose $\{g_n:\partial \mathbb{H}  \rightarrow \partial
\mathbb{H} \}$ and $\{h_n:\partial \mathbb{H}  \rightarrow
\partial \mathbb{H} \}$  uniformly converge to the self
homeomorphisms $g:\partial \mathbb{H}\rightarrow \partial
\mathbb{H}$ and $h:\partial \mathbb{H}\rightarrow \partial
\mathbb{H}$, respectively. Then the sequence $\{f_n\}$ satisfies
axiom (C1) or (C2), if and only if $\{g_nf_nh_n\}$ satisfies (C1)
or (C2), respectively.
\end{lem}

\begin{proof} The (C1) axiom equivalence follows  from the general properties of
topological groups. For the (C2) equivalence, it is enough to show
the property first for postcomposition by $g_n$ and then for
precomposition by $h_n$. To this end,  suppose $\{f_n\}$ locally
uniformly converges to $a$ on $
\partial \mathbb{H}-\{r\}$. Recall that  the spherical metric on $\partial \mathbb{H}$ is
denoted  $\rho(\cdot,\cdot)$. Let $C$ be a compact set in
$\partial \mathbb{H}-\{r\}$ and  let $\epsilon>0$. Choose $\delta
>0$ such that $g$ maps the $\rho$-ball of radius $\delta$ centered
at $a$ into the $\rho$-ball of radius $\epsilon/2$ around $g(a)$.
Choose $n$ sufficiently large so that
$\rho(g_n(x),g(x))<\epsilon/2$, for all $x \in \partial \mathbb
H$, and that $f_n(C)$ lands in the $\rho$-ball of radius $\delta$
centered at $a$. Clearly $g_nf_n(C)$ is in the $\rho$-ball of
radius $\epsilon$ about $g(a)$ for large enough $n$. Thus
$\{g_nf_n\}$ locally uniformly converges to $g(a)$ in $\partial
\mathbb{H} -\{r\}$. The converse follows from group properties. We
have proven the lemma for postcomposition.

To prove the lemma for precomposition, write $f_nh_n$ as
$(h_n^{-1}f_n^{-1})^{-1}$ and use remark (\ref{rem:inverses for
convergence property}).
\end{proof}

\begin{lem} \label{lem:normalization2} Let $\{f_n :\partial \mathbb{H} \rightarrow \partial \mathbb{H}\}$
be a sequence of $K$-quasiconformal mappings.  There exist $g_n
:\partial \mathbb{H} \rightarrow \partial \mathbb{H}$ in
$\text{Conf}(\partial \mathbb{H})$, where
$g_n(f_n(\infty))=\infty$ and $g_nf_n$ is $K$-quasiconformal.
Moreover, a subsequence of $\{f_n\}$ satisfies (C1) or (C2) if and
only if a subsequence of $\{g_nf_n\}$ satisfies (C1) or (C2),
respectively.
\end{lem}

\begin{proof}
Since $\mathcal{K}$ is compact  and acts transitively on $\partial
\mathbb{H}$, we can find $g_n \in \mathcal{K}$ so that $g_n f_n
(\infty)=\infty$, and $g_n$  uniformly converges to $g \in
\mathcal{K}$. Hence the first part of the lemma follows because
the elements of $\mathcal{K}$ are conformal. The last statement
follows from Lemma \ref{lem:normalization1}.
\end{proof}

We omit the   proof of the following lemma, since it closely
follows that of the real case given in \cite{GM}.

\begin{lem}\label{lem:characterize}  Let $\mathcal F \subset \text{Homeo}(\partial \mathbb{H})$ be a family
of self homeomorphisms. Furthermore, assume that $\mathcal F$ is
closed under  post and precomposition by  $\text{Conf}(\partial
\mathbb{H})$.  If $\mathcal F$ has the convergence property then
$\mathcal F$ is a uniformly quasiconformal family with respect to
the conformal structure on $\partial\mathbb H$.
\end{lem}

Using the above lemmas, we have

\begin{thm}\label{thm:convergence prop}  Let $\mathbb{H}$ be a rank one symmetric space of
 noncompact type which is not the real hyperbolic plane.
Then the  group of conformal homeomorphisms of
$\partial\mathbb{H}$ is a convergence group. In fact, any family
$\mathcal F$ of $K$-quasiconformal homeomorphisms of
$\partial\mathbb{H}$ has the convergence property.
\end{thm}

\begin{proof}

To show that the conformal group is a convergence group, set
$\partial \mathbb{H} =\mathcal{N} \cup \{ \infty \}$ and consider
an infinite sequence in $\text{Conf}(\partial \mathbb{H})$. Using
Lemma \ref{lem:normalization2}, it is enough to show that a
sequence $\{g_n \}$ which fixes $\infty$ has an infinite
subsequence which satisfies (C1) or (C2). As remarked in the
basics section, an element of $\text{Conf}(\partial \mathbb{H})$
which fixes $\infty$ acts as a similarity with respect to  the
(Euclidean in the real case) Carnot-Caratheodory metric, $d_c$, on
$\mathcal{N}$. Next, by possibly passing to a subsequence and
inverses, we may assume that the $\{g_n\}$ are  distance
nondecreasing on $\mathcal N$. We continue to call the subsequence
$\{g_n\}$.  According to the Iwasawa decomposition, $g_n \in
\mathcal K \mathcal A \mathcal N$; more precisely, $g_n
=R_n\Lambda_n T_n$, where $R_n \in \mathcal{K}$ fixes $\infty$,
$\Lambda_n$ is a dilation with expansion factor $\lambda_n \geq
1$, and $T_n \in \mathcal{N}$. Since  $\mathcal K$ is compact,
again by passing to a subsequence, we may assume $\{R_n\}$
uniformly converges to $R \in \mathcal{K}$ and hence by Lemma
\ref{lem:normalization1}, it is enough to consider the sequence
$g_n=\Lambda_n T_n$.

 If $\{T_n\}$  has a convergent subsequence then using Lemma
\ref{lem:normalization1}, it's enough to consider the sequence of
dilations $\{ \Lambda_n\}$. But such a sequence clearly has the
convergence property.

Otherwise, if $\{T_n\}$ does not have a convergent subsequence
then $\{T_n\}$ locally uniformly converges to $\infty$ in
$\partial \mathbb{H}-\{\infty\}$. Since $\lambda_n \geq 1$, the
sequence $\{\Lambda_n T_n\}$ shares the same convergence property.

In order to show that any uniformly quasiconformal family
$\mathcal{F}$ has the  convergence property we proceed as follows.
Set $\partial \mathbb{H} =\mathcal{N} \cup \{\infty\}$ and
consider an infinite sequence $\{f_n\} \subset \mathcal F$. Using
Lemma \ref{lem:normalization2} , it is enough to show that a
sequence $\{f_n \}$ of $K$-quasiconformal homeomorphisms of
$\partial \mathbb{H}$ which fix $\infty$ have an infinite
subsequence which satisfies (C1) or (C2). Post compose each $f_n$
by an element $T_n \in \mathcal N$ which takes $f_n(e)$ to $e$.
Next, fix $p \in
\partial \mathbb H -\{e,\infty\}$ and choose an element
$\Lambda_n \in \mathcal A$ taking $T_nf_n(p)$ to the sphere of
radius one centered at the origin.

Set $g_n =\Lambda_n T_n$ and consider the sequence $\{g_n f_n\}$.
Observe that the  $g_nf_n$  fix $e$ and $\infty$ as well as taking
$p$ to the unit sphere centered at the origin. Now by Theorem
(4.8) of \cite{HK}, $\{g_n f_n\}$ being $K$-quasiconformal,
implies $\{g_n f_n\}$ is $K'$-quasisymmetric on $\mathcal{N}$,
where $K'$ only depends on $K$. Therefore, by Theorem (10.26) of
\cite{H}, $\{g_n f_n\}$ has a subsequence (we continue to use the
same notation for the subsequence), which converges locally
uniformly on $\mathcal{N}$ to a map, say $h$. By Exercise (10.29)
of \cite{H}, $h$ is either a $K'$-quasisymmetric embedding or is
constant on $\mathcal{N}$. Since $g_nf_n(p)$ is in the sphere of
radius one, $h$ is not constant. We need to show that
$h:\mathcal{N}\rightarrow \mathcal{N}$ is surjective.

To this end, consider $\mathcal{N}^{\prime}=\partial
\mathbb{H}-\{e\}$. Again $\{g_nf_n\}$ being $K$-quasiconformal
implies $\{g_nf_n\}$ is $K'$-quasisymmetric on
$\mathcal{N}^{\prime}$. The above argument implies $\{g_nf_n\}$
locally uniformly converges to an embedding $h^{\prime}$ of
$\mathcal{N}^{\prime}$. Clearly $h=h^{\prime}$ on $\partial
\mathbb{H}-\{e,\infty\}$ and therefore $h$ is surjective. Thus,
the sequence $\{g_n f_n\}$ uniformly converges to the self
homeomorphism $h$ of $\partial \mathbb{H}$.

Writing $f_n$ as $g_n^{-1}(g_nf_n)$ and appealing to Lemma
\ref{lem:normalization1}, we see that $\{f_n\}$ satisfies (C1) or
(C2) if and only if $\{g_n^{-1}\}$ does. But the elements of
$\{g_n^{-1}\}$ are in $\text{Conf}(\partial \mathbb{H})$ which we
have shown to have the convergence property.
\end{proof}

\begin{thm} \label{thm:maximal} Let $\mathbb{H}_{\mathbb C}$ denote
 complex hyperbolic space of complex dimension greater than one.
 Then the group of conformal homeomorphism of its boundary, $\text{Conf}(\partial \mathbb{H}_{\mathbb C})$, is a maximal convergence group.
 That is,

\begin{enumerate}
\item  $\text{Conf}(\partial \mathbb{H}_{\mathbb C})$ is a
convergence group, and

\item If $G \supseteq \text{Conf}(\partial \mathbb{H}_{\mathbb
C})$ is a convergence group acting on $\partial
\mathbb{H}_{\mathbb C}$, then $G=\text{Conf}(\partial
\mathbb{H}_{\mathbb C})$.
\end{enumerate}
\end{thm}
\begin{proof}

The first item was proven in Theorem \ref{thm:convergence prop}.
For the maximality, suppose $f \in G$ is not conformal. Then by
Lemma 18.5 of \cite{P4}, $f$ is not globally 1-quasiconformal
using the metric definition. That is, the dilatation $H_f(x)>1$ on
a set of positive measure.

The map $f$ is almost everywhere Pansu differentiable (see
$Th\acute{e}or\grave{e}me$ (5) of \cite{P4}). Hence there is a
point of differentiability which we may assume to be $e \in
\partial \mathbb H=\mathcal{N} \cup \infty$ satisfying $H_f(e)>1$.
Using proposition (7) of \cite{KR} we have,

\begin{equation}\label{eq:lambda inequality}
\lambda_{1}^{2n+2}(e) \leq K^{n+1} J_f(e).
\end{equation}

Here $\lambda_1(e)$ is the length of the large axis in the image
of the unit sphere in the contact subspace at $e$ under the
ordinary  derivative of $f$. Also, $ J_f(e)$ is the determinant of
the Pansu derivative at $e$.

Postcomposing by an element of $\text{Conf}(\partial
\mathbb{H}_{\mathbb C})$, we may assume that  $f(e)=e$ and
$J_f(e)$ is one. Moreover, since the stabilizer of $e$ in
$\text{Conf}(\partial \mathbb{H}_{\mathbb C})$ acts transitively
on directions in the contact subspace at $e$, we may assume that
$\lambda_{1}(e) $ is an eigenvalue of $f$ with eigenvector in the
contact subspace. Note that by the chain rule all iterates of $f$
have $e$ as a point of Pansu differentiability. Thus by inequality
(\ref{eq:lambda inequality}), iterating $f$ gives rise to
quasiconformal maps with arbitrarily large dilatation in $G$,
which is a contradiction. Therefore, $G=\text{Conf}(\partial
\mathbb{H}_{\mathbb C})$.
\end{proof}

For real hyperbolic space of dimension greater than two, the fact
that the M\"obius group is a maximal convergence group acting on
$\partial \mathbb{H}$ was proven in Gehring-Martin \cite{GM}. In
the case of  Quaternionic or Cayley Hyperbolic space, all
quasiconformal self homeomorphisms of $\partial
\mathbb{H}_\mathbb{H}^{4n}$ ($n \geq 2$), and $\partial
\mathbb{H}_{\mathbb{C}a}^{16}$ are 1-quasiconformal, and hence
conformal. This follows from (11.2) and (11.5) of \cite{P3}.
Putting these results together with Theorem (\ref{thm:maximal}) we
have,

\begin{cor} \label{cor:maximal}Let $ \mathbb H$  be  a rank one symmetric
space of noncompact type  different from the real hyperbolic
plane. Then the group of conformal homeomorphisms,
$\text{Conf}(\partial \mathbb{H})$, is a maximal convergence
group.
\end{cor}

\section{Some Consequences of Maximality}
We remind the reader that $\mathbb{H}$  denotes a rank one
symmetric space of noncompact type except for the hyperbolic
plane. The following fact is an immediate application of Corollary
\ref{cor:maximal}.

\begin{cor} \label{cor:qcmaximal} The group of
conformal homeomorphisms, $\text{Conf}(\partial \mathbb{H})$, is a
maximal uniformly quasiconformal group.
\end{cor}

Let $\text{QC}(\partial \mathbb{H})$ denote the group of
quasiconformal homeomorphisms of the boundary. A {\it full
quasi-isometry} of $\mathbb{H}$ is a quasi-isometry $f:\mathbb{H}
\rightarrow \mathbb{H}$ whose domain and image is $\mathbb{H}$;
two such are equivalent if the Hausdorff distance between their
graphs in $\mathbb{H}\times \mathbb{H}$ is finite.  Let
$\text{QI}(\mathbb{H})$ denote the group of equivalence classes of
full quasi-isometries of $\mathbb{H}$. Note that
$\text{Isom}(\mathbb{H})$ naturally embeds in
$\text{QI}(\mathbb{H})$.  A family ${\mathcal{F}} \subset
\text{QI}(\mathbb{H})$ is said to be a {\it uniformly
quasi-isometric} family if there exists a uniform bound on the
Lipschitz constant for all elements of ${\mathcal{F}}$. It is well
known that a quasi-isometry class $[f]$ induces a quasiconformal
homeomorphism, $f|_{\partial H}$, on the boundary. Moreover, the
quasiconformality constant is dependent solely on the Lipschitz
constant of the quasi-isometry (see  \cite{P3}). In fact, the map
$[f] \mapsto f|_{\partial H}$, is a group isomorphism between $
\text{QI}(\mathbb{H})$ and $\text{QC}(\partial \mathbb{H})$, which
identifies $\text{Isom}(\mathbb{H})$ with $\text{Conf}(\partial
\mathbb{H})$. Moreover, $\mathbb{H}^{4n}_\mathbb{H}$, for $n \geq
2$, and $\mathbb{H}^{16}_{\mathbb{C}a}$ are $\text{QI}$-rigid;
that is,
$\text{Isom}(\mathbb{H}^{4n}_\mathbb{H})=\text{QI}(\mathbb{H}^{4n}_\mathbb{H})$
and
$\text{Isom}(\mathbb{H}^{16}_{\mathbb{C}a})=\text{QI}(\mathbb{H}^{16}_{\mathbb{C}a})$.
In these rigid cases, maximality of $\text{Isom}(\mathbb{H})$ as a
uniform quasi-isometry group is a triviality. A more interesting
phenomenon is that of the maximality of $\text{Isom}(\mathbb{H})$
in the remaining cases of $\mathbb{H}_\mathbb{R}^n$, for $n \geq
3$, and $\mathbb{H}_\mathbb{C}^{2n}$, for $n \geq 2$ . This fact,
which has been observed in Gromov and Pansu, \cite{GP}, may be
thought of as a simple consequence of Corollary
\ref{cor:qcmaximal}.

\begin{cor}\label{cor:psl(2,R)maximalqi}
The group of isometries, $\text{Isom}(\mathbb{H})$, is a maximal
uniformly quasi-isometric group.
\end{cor}

Another elementary consequence of Corollary  \ref{cor:maximal} is
the following,

\begin{cor}\label{prop:isom(H)maximalconv}
 Let $G$ be a convergence group acting on $\mathbb{H}\cup
\partial\mathbb{H}$. Suppose that each element of $G$ is topologically conjugate
to an element of $\text{Isom}(\mathbb{H})$. If
$\text{Isom}(\mathbb{H})\leq G$, then $G=\text{Isom}(\mathbb{H})$.
The conjugating homeomorphism need not be the same for all
elements of $G$.
\end{cor}

\begin{proof} Consider
the map $\phi: G \rightarrow Homeo(\partial\mathbb{H})$, given by
$f \mapsto f|_{\partial\mathbb{H}}$. Note that
$\text{Image}(\phi)$ is a convergence group acting on
$\partial\mathbb{H}$ and contains
$\text{Conf}(\partial\mathbb{H})$. Since
$\text{Conf}(\partial\mathbb{H})$ is a maximal convergence group
(Corollary \ref{cor:maximal}) we have $\text{Image}(\phi)
=\text{Conf}(\partial\mathbb{H})$.

Now take $f \in Ker(\phi)$. By our assumption, there exists a
homeomorphism $h: \mathbb{H}\cup \partial\mathbb{H} \rightarrow
\mathbb{H}\cup \partial\mathbb{H}$, such that $h^{-1} \circ f
\circ h \in \text{Isom}(\mathbb{H})$.  Note that $\phi(h^{-1}
\circ f \circ h) = h^{-1}|_{\partial\mathbb{H}}\circ
f|_{\partial\mathbb{H}} \circ h|_{\partial\mathbb{H}}
=h^{-1}|_{\partial\mathbb{H}}\circ Id|_{\partial\mathbb{H}} \circ
h|_{\partial\mathbb{H}}= Id|_{\partial\mathbb{H}}$. Since $\phi$
restricted to $\text{Isom}(\mathbb{H})$ is a monomorphism, it must
be that $h^{-1} \circ f \circ h=Id_{\mathbb{H}}$ and hence $f =
Id_{\mathbb{H}}$;  that is, $\phi : G \rightarrow Homeo(\partial
\mathbb{H})$ is a monomorphism. Since
$\text{Image}(\phi)=\phi(\text{Isom}(\mathbb{H}))=
\text{Conf}(\partial\mathbb{H})$, we may conclude
$G=\text{Isom}(\mathbb{H})$.
\end{proof}

\end{document}